\documentclass[11pt]{amsart}

\usepackage[toc,page]{appendix} 
\usepackage{amsfonts} 
\usepackage{amssymb}
\usepackage[final]{graphicx, graphics,stmaryrd} 
\input xy
\xyoption{all}
\usepackage[margin=3cm]{geometry} 
\usepackage{amscd}
\usepackage[all]{xy}

\usepackage[breaklinks,bookmarksopen,bookmarksnumbered]{hyperref}

\usepackage{tensor}
\usepackage[english]{babel}
\usepackage[utf8]{inputenc}
\usepackage{amsthm}
\usepackage{graphics}
\usepackage{amsmath}
\numberwithin{equation}{section} 
\usepackage{amstext}
\usepackage{engrec}
\usepackage{rotating}
\usepackage[safe,extra]{tipa}
\usepackage{multirow}
\usepackage{enumitem} 
\usepackage{bm}
\usepackage{stackrel} 
\usepackage{tikz-cd} 

\usepackage{algorithmic,algorithm,mathtools} 

\newtheorem{teo}{Theorem}[section]

\newtheorem{prop}[teo]{Proposition}
\newtheorem{lemma}[teo]{Lemma}

\theoremstyle{definition}
\newtheorem{defin}[teo]{Definition}
\newtheorem{rmk}[teo]{Remark}

\newcommand{\mc}{\mathcal}

\def\R{\mathbb{R}}

\def\Z{\mathbb Z}

\def\P{\mathcal P}
\def\cc{\mathcal C}

\def\cc{\mathcal{C}}

\def\H{\mc{H}}

\newcommand{\m}{\mbox}

\newcommand{\cor}{\textit}

\newcommand{\fine}{\qed\newline}

\newcommand{\beginCD}{\begin{equation*}\begin{CD}}
\newcommand{\enCD}{\end{CD}\end{equation*}}

\DeclareMathOperator{\id}{id}

\DeclareMathOperator{\argmin}{argmin}

\newcommand{\al}{\alpha}

\DeclareMathOperator{\adm}{Adm}
\DeclareMathOperator{\Opt}{Opt}

\def\mmu{\bm{\mu}}

\DeclareMathOperator{\Geod}{Geod}
\DeclareMathOperator{\supp}{supp}
\DeclareMathOperator{\iter}{iter}

\def\de{\partial}

\newcommand{\norm}[1]{\left\lVert#1\right\rVert}

\title{Cortically based optimal transport}
\author{Mattia Galeotti, Giovanna Citti, Alessandro Sarti}

\begin{document}
\maketitle

\begin{abstract}
We introduce a model for image morphing in the primary visual cortex V1 to perform completion of missing images in time. We model the output of simple cells through a family of Gabor filters and the propagation of the neural signal accordingly to the functional geometry induced by horizontal connectivity. 
Then we model the deformation between two images as a path relying two different outputs. This path is obtained by optimal transport considering the Wasserstein distance geodesics
associated to some probability measures naturally induced by 
the outputs on V1. The frame of Gabor filters allows to project back the output path, therefore obtaining an associated image stimulus deformation.
We perform a numerical implementation of our cortical model,
assessing its ability in reconstructing rigidi motions
of simple shapes.
\end{abstract}

\section{Introduction}

The functional geometry of the visual cortex is a widely studied subject. 
It is known that cells of the primary visual cortex V1 are sensitive to specific features 
of the visual stimulus, like position, orientation, scale, colour, curvature, velocity and many others \cite{hub88}. 
In the seventies the neurophysiologists Hubel and Wiesel discovered the modular organisation of the primary visual cortex \cite{hubwie77}, 
meaning that cells are spatially organized in such a way that for every point $(x,y)$ of the retinal plane 
there is an entire set of cells, each one sensitive to a particular instance of the considered feature. 
This organisation corresponds to the so-called hypercolumnar structure. 
Hypercolumns of cells are then connected to each other by means of the horizontal connectivity, 
allowing cells of the same kind but sensitive to different points $(x,y)$ of the stimulus  to communicate. 
Hypercolumnar organization and neural connectivity between hypercolumns constitute the 
functional architecture of the visual cortex, 
that is the cortical structure underlying the low level processing of the visual stimulus.
The mathematical modelling of the functional architecture of the visual cortex in terms of 
differential geometry was introduced in the seminal works of Hoffmann \cite{hoff66, hoff89}, 
who proposed to model the hypercolumnar organization in terms of a fiber bundle. 
Many of such results dealing with differential geometry were given 
a unified framework under the new name of neurogeometry.\newline

 Petitot and Tondut \cite{petiton99}, related the contact geometry introduced by Hoffmann 
 with the geometry of illusory contours of Kanizsa \cite{kani80}. 
 The problem of completion of occluded object was afforded 
 by computing geodesic curves in the contact structure.
 
Then, in \cite{cisar06} Citti and Sarti showed how the functional architecture 
could be described in terms of Lie groups structures. 
In particular, as proved by Daugman \cite{daug85} the receptive profiles of simple cells 
can be modelled as a Gabor filter. Since these filters can be obtained via rotation and translation 
from a fixed one, the functional architecture of the whole family of simple cells 
have been described as the Euclidean motion group $SE(2)$ \cite{cisar06}. 
In presence of a visual stimulus $I: \R^2 \to  [0, 1]$ on the retinal plane, 
the action of the whole family of simple cells is obtained by convolving the function 
$I$ with the bank of Gabor filters. 
The output of the cells action will be a function $\mu: \R^2 \times S^1 \to \R$. 
The horizontal connectivity is strongly anisotropic and it is modelled via a sub-Riemannian metric. 
Since it is very common that part of the visual input is occluded,
the action of horizontal connectivity allows to complete the missing part by means of propagation 
in such a space, 
under the action of advection diffusion differential operators of Fokker-Planck type.
Visual completion problems are then solved via geodesics or minimal surfaces. 
This approach was extended to scale in \cite{scp08}, to space-time in \cite{bccs14} and to frequency in \cite{bsc19}. 
In \cite{fahs16} and \cite{scp19} the lifting has been extended to heterogeneous features defined in different groups.
For an extended review on neurogeometry see \cite{cisar14}.\newline

In this paper we aim to reconsider the problem of  completion of missing stimulus in time by means of morphing of one lifted cortical image 
in a different one in terms of optimal transport of a probability distribution in the functional geometry of the cortex. 
Two images can represent the same object at two different intervals of times, 
and different algorithms have been proposed to perform completion of missing images between the two.
We recall the results of \cite{wsi04} and the model proposed by \cite{zyht07} on the image plane as a geodesic in the Wasserstein space. 

We propose a cortical version of this phenomena, using geodesics in the cortical space endowed 
with the Wasserstein distance. We work in the manifold $M = \R^2 \times S^1\times \R^+$ 
of cells sensible to position orientation and scale (see also \cite{scp08}),
and we develop a model that treats well shape rotations.
We point out that the sub-Riemannian metric in $\R^2\times S^1$
 is important to keep together the shape of the object along the rotational movements. 
 In fact the distance function in this metric approximates the statistical correlations of boundaries in natural images \cite{scs10}.

In order to obtain a metric deformation, we consider the positive
and negative part of  the output $\mu$ at a given time, normalized as two probability
measures $\mu^+,\mu^-$ over $M$. Following
the papers \cite{amgi08, amgi13} of Ambrosio and Gigli,
we consider the space $\P_2(X)$ of probabilities
with finite $2$-momentum over $X=\R^2\times S^1\times \R$,
endowed with the associated Wasserstein distance; on this
space we can find, for any pair of inputs $I_0,I_1$ taken at different times, a unique
constant speed geodesic relying their
associated output measures.

In particular, for any regular measure $\mu\in \P_2(X)$, we obtain
Theorem \ref{teo_oneT}, a generalization of previously known results,
and it
assures that for any measure $\nu\in \P_2(X)$ there exists
a unique transport map $T$ in the sense of Monge's
formulation of optimal transport. Using this transport map,
it is possible to give an explicit description of the constant speed
geodesics in~$\P_2(X)$; in the case of the measures induced by the output functions,
this is done in Remark \ref{rmk_et} and equation (\ref{path}).
Using the frame properties of the family of Gabor filters,
we reconstruct with equation (\ref{eq_imp0}) a path of images $I_t$
relying the initial and final input $I_0,I_1$.

In the two final sections, we develop a numerical implementation
of our model. In fact, the frame generated by the 
odd Gabor mother function is not invertible in a discrete setting,
and in order to make our model workable, we use a Wavelet Gabor Pyramid
generated by an odd and an even Gabor function.
The outputs obtained via these frames are transported
following the same approach detailed above. As shown in Section \ref{impl},
our model allows the deformation of simple shapes
through rotation and translation, preserving the basic structures
in the treated pictures. As we
discuss, the same preservation is not attained
by the numerical 2-dimensional `classical' implementation 
of optimal transport, meaning the implementation of the optimal
transport between the inputs $I_0,I_1$ seen as measures on $\R^2$
and taking the square Euclidean distance as the cost function.\newline

In Section \ref{retina} we introduce the operation
associating to any input an output function via the convolution
with the family of Gabor filters. In Section \ref{optimal}
the classical problem of optimal transport is introduced,
with the techniques necessary for a general solution.
In Section \ref{sec_geo} we describe the constant speed geodesics
in~$\P_2(X)$. In Section \ref{recon}
we treat the properties of continuous frames such as the one
of Gabor functions; they are fundamental in reconstructing
a path of input images from the deformation path of output functions.
Section \ref{defoutput} states that the measure deformation results
are valid in our case. Finally, in Section \ref{const}
we find a constraining condition implying that the path
of output functions $\mu_t$ is naturally induced by a path
of input images $I_t$. In Section \ref{discmod} we introduce the
discrete setting and the Wave Gabor Pyramid. Finally in Section \ref{impl}
we discuss the numerical implementation and its results.

\section{From the retina to the output space}\label{retina}
We consider an input function $I\colon \R^2\to [0,1]$ in $L^2(\R^2)$. This functions
models an input received in the retina plan
and it induces an ouput function on the cortex. 
In order to introduce the output function, we recall 
the odd series of Gabor filters.
We call Gabor mother filter the function 
$$\psi_{0,0,0,1}\colon (\tilde x,\tilde y)\in \R^2\mapsto e^{-(\tilde x^2+\tilde y^2)}\cdot\sin(2\tilde y),$$
Moreover, we consider the roto-dilation defined by
$$A_{\theta,\sigma}:=\sigma\cdot R_\theta=\sigma\cdot\left(\begin{array}{cc} \cos\theta &-\sin\theta\\
\sin\theta &\cos\theta\\
\end{array}\right),$$
for any $\theta \in S^1$ and $\sigma \in \R^+$.
We also consider the application 
$$A_{x,y,\theta,\sigma}:(\tilde x,\tilde y)\mapsto A_{\theta,\sigma}(\tilde x,\tilde y)+(x,y).$$
All this allows the definition of a family of Gabor filters
$$\psi_{x,y,\theta,\sigma}(\tilde x,\tilde y):=\frac{1}{\sigma^{3\slash2}}\cdot\psi_{0,0,0,1}(A_{\theta,\sigma}^{-1}(\tilde x-x,\tilde y-y))=
\frac{1}{\sigma^{3\slash2}}\cdot\psi_{0,0,0,1}(A_{x,y,\theta,\sigma}^{-1}(\tilde x,\tilde y)).$$

In what follows we will denote the filter $\psi_{0,0,0,1}$ by $\psi_0$ when there is no risk of confusion.
We consider the variety $M:= \R^2\times S^1\times \R^+$ with
its natural (Lebesgue) measure $dk$, this is the output space
where we build the $\mu$ function induced by the input $I$.
For any point $k=(x,y,\theta,\sigma)\in M$,
we denote by $\psi^k$ the Gabor filter $\psi_{x,y,\theta,\sigma}$.\newline

\begin{defin}\label{out}
Consider a point $k=(x,y,\theta,\sigma)$
on $M$, then the output function of a cell
in response to the visual input $I$ is
\begin{align*}
\mu(x,y,\theta,\sigma)=\mu(k)&:=\langle I,\psi^k\rangle\\
&=\int_{\R^2}I(\tilde x,\tilde y)\psi_{x,y,\theta,\sigma}(\tilde x,\tilde y)d\tilde xd\tilde y.\\
\end{align*}
\end{defin}

\begin{rmk}\label{rmk_met}
We work in the case $M=\R^2\times S^1\times \R^+$
and we put on $M$ the Riemannian structure which endows the neurogeometry
of the cortex. In particular the $\R^2$ factor is the retinal plan with the natural projection
$M\to \R^2$, $\theta\in S^1$ is an angular parameter
that encodes the border orientation in the processes of
border recognition, while $\sigma\in \R^+$
is a scale parameter in the same process.

In order to define this metric $g$, consider
the following four vector fields, in every point
they span the $M$ tangent bundle,
\begin{align*}
X_1&=\cos\theta \cdot\partial_x+\sin\theta\cdot \partial_y,\\
X_2 &=\partial_\theta,\\
X_3 &=-\sin\theta\cdot\partial_x+\cos\theta\cdot\partial_y,\\
X_4 &= \partial_\sigma.
\end{align*}
Let $g$ be the metric such that at any point
$ \frac{1}{\sqrt\sigma}X_1,\ \frac{1}{\sqrt\sigma}X_2,\  \sqrt\sigma X_3,\ \sqrt\sigma X_4$ is an orthonormal system.
That is at every point the metric $g$ is represented by the matrix
$$
\tilde g=\left(\begin{array}{cccc}
\sigma &0 &0 &0\\
0 &\sigma &0&0\\
0 &0 &\frac{1}{\sigma} &0\\
0& 0& 0& \frac{1}{\sigma}\\
\end{array}\right)
$$
For $\sigma$ that tends to $0$, $g$ tends
to the sub-riemannian structure treated in \cite{cisar06} (module rescaling), for $\sigma\to + \infty$ 
it approaches a hyperbolic metric.

In this setting the measure $\mu\cdot dk=\langle I,\psi^k\rangle \cdot dk$ on $M$, is
$\langle I,\psi^k\rangle$ times the measure induced by the metric $g$.
\end{rmk}

As proved in Appendix \ref{appa}, the integral 
$$\int_{M}\psi^k(\tilde x,\tilde y)dk$$
has finite value $0$, independently of the pair $(\tilde x,\tilde y)$.

\begin{rmk}\label{rmk_mass}
The mass of the output function $\mu$ on $M$, is proportional  to the mass of~$I$
and therefore it is null. Indeed,
\begin{align*}
\int_{M}\mu(k)dk&=\int_{M}dk\int_{\R^2}I(\tilde x,\tilde y)\psi^k(\tilde x,\tilde y)d\tilde xd\tilde y\\
&=\int_{\R^2}I(\tilde x,\tilde y)d\tilde xd\tilde y\int_{M}\psi^k(\tilde x,\tilde y)dk=0
\end{align*}
by the Fubini's Theorem. \newline
\end{rmk}

\section{The optimal transport problem}\label{optimal}

We recall the classical Kantorovich's formulation of optimal transport. Our
main reference is the Ambrosio-Gigli guide \cite{amgi13}.
For any Polish space $X$ (\cor{i.e.}~a complete and separable metric space)
we denote by $\P(X)$ the set of Borel probability measures on $X$
and by $\mc B(X)$ the set of Borel sets on $X$.
Consider two Polish spaces $X,Y$,
if $\mu\in \P(X)$ and $T\colon X\to Y$
is a Borel map, then we denote by $T_\#\mu\in \P(Y)$
the pushforward of $\mu$ through $T$, defined by
$$T_\#\mu(E)=\mu(T^{-1}E)\ \ \ \forall E\in \mc B(Y).$$

Consider the natural product $X\times Y$ and its associated projections $\pi^X,\pi^Y$.
Let $c\colon X\to Y$ be a Borel map called \cor{cost} function, and consider two measures $\mu\in\mc{P}(X)$
and $\nu\in\P(Y)$.
\begin{defin}
An \cor{admissible transport plan} between $\mu$ and $\nu$ is a measure $\gamma\in\mc{P}(X\times Y)$
such that $\pi^X_\#\gamma=\mu$ and $\pi^Y_\#\gamma=\nu$, or equivalently
\begin{align*}
\gamma(A\times Y)&=\mu(A)\ \ \forall A\in\mc{B}(X)\\
\gamma(X\times B)&=\nu(B)\ \ \forall B\in\mc{B}(Y).
\end{align*}
We denote the set of admissible transport plans between $\mu$ and $\nu$ by
$\adm(\mu,\nu)$.
\end{defin}

We want to minimize the integral 
$$\int_{X\times Y}c(x,y)d\gamma(x,y)$$
for all the admissible transport plans between $\mu$ and $\nu$. We say that $\gamma$ is induced
by a transport map if there exists a Borel map $T\colon X\to Y$ such that $\gamma=(\id\times T)_\#\mu$,
in that case
$$\int_{X\times Y}c(x,y)d\gamma=\int_Xc(x,T(x))d\mu.$$
An optimal transport plan is a transport plan $\gamma$ that realizes the infimum above.
It is known that such a minimizer exists under very general conditions. We denote
the set of optimal plans by $\Opt(\mu,\nu)$.
\begin{prop}[see {\cite[Theorem 4.1]{villa08}}]
Consider $\mu\in\P(X)$ and $\nu\in\P(Y)$. 
If the cost function $c$ is lower semicontinuous and bounded from below, then there exists
an optimal plan $\gamma$ for the functional
$$\gamma\mapsto \int_{X\times Y}c(x,y)d\gamma(x,y),$$
among all $\gamma\in \adm(\mu,\nu)$.
\end{prop}

We are interested in the cases where an optimal plan is induced by a transport map $T$.
We state  \cite[Lemma 1.20]{amgi13}, referring to Ambrosio-Gigli paper for a proof.
\begin{lemma}
Consider $\gamma\in\adm(\mu,\nu)$. Then $\gamma$ is induced by a map
if and only if $\gamma$ is concentrated in a measurable set $\Gamma\subset X\times Y$
such that for $\mu$-a.e.~$x$ there exists only one $y=T(x)$ in $\Gamma\cap (\{x\}\times Y)$.
In this case, $T(x)$ induces $\gamma$.\newline
\end{lemma}

In order to introduce the notion of
$c$-concavity for some cost function $c$, 
and in order to show the existence and uniqueness of an optimal 
transport plan $\gamma$ induced by a transport map $T$, we 
give some definitions following \cite[Chap.1]{amgi13}.

\begin{defin}[Superdifferential]\label{def_sd}
Consider $M$ a riemannian manifold and any function $\varphi\colon M\to \R$, 
we define its superdifferential at any point $x\in M$,
$$\de^+\varphi(x):=\left\{dh(x)\in T_x^*M:\ h\in\cc^1(M,\R),\ \varphi-h\m{ attains a local maximum at }x\right\}.$$
When there is no risk of confusion, we will denote by $\de^+\varphi$ the
associated subspace of the total space $T^*M$.
The subdifferential $\de^-\varphi$ is defined analogously as
the set of differentials $dh(p)$ where $\varphi-h$ attains a local minimum of $\varphi-h$.
\end{defin}
\begin{rmk}
Equivalently, $\de^+\varphi(x)$ is the set of vectors
$v\in T_xM$ such that
$$\varphi(z)-\varphi(x)\leq \langle v,\exp^{-1}_x(z)\rangle+o(d(x,z)).$$
The same for $\de^-\varphi(x)$ with inversed inequality.
With this definition, $\de^+\varphi$ and $\de^-\varphi$ are subspaces of $TM$.
\end{rmk}
It is well known that if $\varphi$ is differentiable at $x\in M$, its superdifferential
and subdifferential at $x$ coincide and the contain only the $\varphi$ gradient,
$$\de^+\varphi(x)=\de^-\varphi(x)=\{\nabla\varphi(x)\}.$$\newline

Consider two Polish spaces $X,Y$, and a cost function $c\colon X\times Y\to \R$.
\begin{defin}[$c$-transforms]
Consider a function $\varphi\colon X\to \R\cup\{\pm\infty\}$, its $c_+$-transform $\varphi^{c_+}\colon Y\to \R\cup\{\pm\infty\}$
is defined as
$$\varphi^{c_+}(y):=\inf_{x\in X}c(x,y)-\varphi(x).$$
Analogously for any $\psi\colon Y\to \R\cup\{\pm\infty\}$, we can define its $c_+$-transform $\psi^{c_+}\colon X\to \R\cup\{\pm\infty\}$.

The $c_-$-transform of $\varphi$ is $\varphi^{c_-}\colon Y\to \R\cup\{\pm\infty\}$ defined as
$$\varphi^{c_-}(y):=\sup_{x\in X}-c(x,y)-\varphi(x).$$
Analogously for the $c_-$-transform $\psi^{c_-}$ of $\psi$.
\end{defin}

\begin{defin}[$c$-concavity]
We say that a function $\varphi\colon X\to \R\cup\{-\infty\}$ is $c$-concave
if there exists $\psi\colon Y\to\R\cup\{-\infty\}$ such that $\varphi=\psi^{c_+}$.
Analogously we have a notion of $c$-convexity.
\end{defin}

\begin{defin}[Semiconcavity]
A function $f\colon U\to \R$ whose domain is a convex subset $U$ of a riemannian manifold $M$,
is semiconcave with constant $K$ if for every geodesic $\gamma\colon [0,1]\to U$ and $t\in [0,1]$,
we have
$$t\cdot f(\gamma_0)+(1-t)\cdot f(\gamma_1)\leq f(\gamma_t)+\frac{1}{2}\cdot t(1-t) K\cdot d^2(\gamma_0,\gamma_1).$$
\end{defin}
With the notation of \cite[Chap.10]{villa08}, the definition above describes a 
semiconcave function with modulus $\omega(t)=K\frac{t^2}{2}$.

\begin{rmk}\label{rmk_V}
Observe that by \cite[Equation (10.14)]{villa08}, if $f$ is
semiconcave, superdifferentiable at $x$
and $q\in \de^+f(x)$, then 
$$f(\exp_x w)-f(x)\leq \langle q, w\rangle -\frac{1}{2}K\norm{w}^2.$$
\end{rmk}

\begin{defin}[$c$-superdifferential]
Consider $\varphi\colon X\to \R\cup\{-\infty\}$ a $c$-concave function, 
then its $c$-superdifferential $\de^{c_+}\varphi\subset X\times Y$ is defined as
$$\de^{c_+}\varphi:=\left\{(x,y):\  \varphi(x)+\varphi^{c_+}(y)=c(x,y)\right\}.$$
We denote by $\de^{c_+}\varphi(x)$ the set of $y\in Y$ such that $(x,y)\in \de^{c_+}\varphi$.
Analogously we can define the $c$-subdifferential $\de^{c_-}\varphi\subset X\times Y$.
\end{defin}

Consider two probability measures $\mu\in \P(X)$ and $\nu\in \P(Y)$.
In the following we will consider a cost function
$c\colon X\times Y\to \R$ such that there exists
two 
functions $a\in L^1(\mu)$, $b\in L^1(\nu)$,
respecting the inequality below.
\begin{equation}\label{condition}
c(x,y)\leq a(x)+b(y).
\end{equation}

\begin{teo}[{\cite[Theorem 1.13]{amgi13}}]\label{teo_supp}
Consider $\mu\in\P(X)$, $\nu\in\P(Y)$ and $c\colon X\times Y\to \R$ a continuous and bounded from below cost function
such that there exist two functions $a\in L^1(\mu)$ and $b\in L^1(\nu)$ verifying
condition (\ref{condition}).
Then there exists a $c$-concave function $\varphi\colon X\to \R$ such that
$\varphi\in L^1(\mu)$, $\varphi^{c_+}\in L^1(\nu)$
and for any optimal plan $\gamma\in \Opt(\mu,\nu)$,
$$\supp(\gamma)\subset \de^{c_+}\varphi.$$
\end{teo}

Consider the manifold $X=\R^2\times S^1\times \R$ with the Lebesgue metric.
Consider a function $\varphi\colon X\to \R$
and the square distance
$c(x,y)=d^2(x,y)\slash2$ as cost function, then
the following lemma and proposition state the link between the superdifferential $\de^+\varphi$
and the $c$-superdifferential $\de^{c_+}\varphi$.

\begin{lemma}\label{lems}
For any point $y\in X$, the function $d^2(-,y)\slash 2$ is uniformly
semiconcave on~$X$.
\end{lemma}

This is proven following the line of reasoning of \cite[Third Appendix]{villa08},
because $X$ is flat and therefore its sectional curvature is everywhere $0$.\newline

\begin{prop}\label{prop_conc}
Consider a $c$-concave function $\varphi\colon X\to \R$, then
it must be 
semiconcave. 
Furthermore, for any $x\in X$, $\exp_x^{-1}(\de^{c_+}\varphi(x))\subset -\de^+\varphi(x)$.
\end{prop}

\proof This is a slight generalization of \cite[Proposition 1.30]{amgi13},
and we develop the same arguments.
As stated by Lemma \ref{lems}, for any $y\in X$ the distance function
$d^2(-,y)\slash 2$ is 
semiconcave,
therefore by Remark \ref{rmk_V} this implies
$$\frac{d^2(z,y)}{2}-\frac{d^2(x,y)}{2}\leq -\langle v,\exp_x^{-1}(z)\rangle +o(d(x,z)),$$
because if $v\in \exp_x^{-1}(y)$, then $-v$ is in the superdifferential of $d^2(-,y)\slash 2$ at $x$.

If we take $d^2(x,y)\slash 2$ as the cost function $c(x,y)$
and consider $y\in \de^{c_+}\varphi(x)$, therefore by definition
$\varphi(z)-c(z,y)\leq \varphi(x)-c(x,y)$ for any $z\in X$. As a consequence
$$\varphi(z)-\varphi(x)\leq \frac{d^2(z,y)}{2}-\frac{d^2(x,y)}{2}\leq \langle -v,\exp_x^{-1}(z)\rangle +o(d(x,z)),$$
that means $-v\in \de^+\varphi(x)$.\fine

\begin{defin}[Regular measure]
We say that a measure $\mu\in \P(X)$ is regular if it vanishes on the set of points of non differentiability of
any semiconcave function $\varphi\colon X\to \R$.
\end{defin}

For any Polish space $(X,d)$, we introduce the space
of probability measures on $X$ with finite $2$-momentum with respect to $d$,
$$\P_2(X):=\left\{\mu\in\P(X):\ \int d^2(x,x_0)d\mu<\infty\ \m{for some, and thus any, }x_0\in X\right\}.$$
We regard $\P_2(X)$ as a metric space with respect to the sup norm.

\begin{teo}\label{teo_oneT}
Consider the riemannian manifold $X$ and a probability measure $\mu\in\P_2(X)$.
If $\mu$ is regular, then for every $\nu\in\P_2(X)$ there exists only one transport plan from $\mu$ to $\nu$
and it is induced by a map $T$. If this is the case, the map $T$ can be written as
$x\mapsto \exp_x(-\nabla\varphi(x))$ for some $c$-concave function $\varphi\colon X\to \R$.
\end{teo}

Observe that this is a generalization of \cite[Theorem 1.33]{amgi13}
to the case of the non-compact riemannian manifold $X$. For another
reference see also \cite{ganmcc96}.

\proof In order to apply Theorem \ref{teo_supp} we are going to verify
that condition (\ref{condition}) is respected.
In particular we want to show that taking $a(x)=d^2(x,x_0)$
and $b(y)=d^2(y,x_0)$ for any $x_0\in M$,
then $a\in L^1(\mu)$, $b\in L^1(\nu)$ and
$c(x,y)=d^2(x,y)\slash 2\leq d^2(x,x_0)+d^2(y,x_0)$.
The inequality is proved by 
\begin{align*}
d^2(x,x_0)+d^2(y,x_0)&\geq\frac{d^2(x,x_0)+d^2(y,x_0)}{2}+|d(x,x_0)d(y,x_0)|\\
&=\frac{1}{2}(d(x,x_0)+d(y,x_0))^2\geq\frac{d^2(x,y)}{2}.
\end{align*}
To have $a\in L^1(\mu)$ (and therefore $b\in L^1(\nu)$) it suffices to have
$\int_Md^2(x,x_0)\mu(x)< \infty$ but this is exactly the definition of $\mu\in \P_2(X)$.

 Thus as a consequence of Theorem \ref{teo_supp},
there exists a $c$-concave function $\varphi$ such that any optimal plan
$\gamma$ is concentrated on $\de^{c_+}\varphi$.
By Proposition~\ref{prop_conc}, $\varphi$
is semiconcave and therefore differentiable $\mu$-a.e.~by $\mu$-regularity.
If $\varphi$ is differentiable at $x$, then $\de^+\varphi(x)$ is the singleton $\{\nabla\varphi(x)\}$,
therefore $\de^{c_+}\varphi(x)$ is empty or equals $\exp_x(-\nabla\varphi(x))$. 
We define $\mu$-a.e.
the function $T(x)=\exp_x(-\nabla\varphi(x))$. As $\supp(\gamma)\subset\de^{c_+}\varphi$,
we must have that $T$ induces $\gamma$, concluding the proof.\fine

\section{Geodesics in ${\mc{P}_2(X)}$}\label{sec_geo}
The introduction of the Wasserstein distance $W_2$ on $\mc P_2$ allows
the definition of geodesics in this space. We show a theorem 
of existence and uniqueness for such a geodesic relying a pair
of measures.

\begin{defin}
Consider a metric space $(X,d)$.
A curve $\gamma\colon [0,1]\to X$ is a constant speed geodesic  if
$$d(\gamma_s,\gamma_t)=|t-s|d(\gamma_0,\gamma_1)\ \ \forall t,s\in [0,1].$$
\end{defin}
We recall that $(X,d) $
is called a geodesic space if
for every $x,y\in X$ there exists a constant speed geodesic connecting them.
We consider the metric space $\Geod(X)$ of constant speed geodesics endowed
with the sup norm.

On $\Geod(X)$ we introduce for any $t\in [0,1]$ the map $e_t\colon \Geod(X)\to X$
such that
$$e_t\colon \gamma\mapsto \gamma_t.$$
Furthermore we define the Wassertein distance associated to $d$ on $\P_2(X)$.

\begin{defin}[Wasserstein distance]
If $\mu,\nu\in\P_2(X)$, then
$$W_2^2(\mu,\nu):=\inf_{\gamma\in\adm(\mu,\nu)}\int d^2(x,y)d\gamma.$$
\end{defin}

The following theorem is proved in \cite{amgi13}.
For any two probability measure $\mu_0,\mu_1$ on a Polish space,
this gives a constant speed geodesic relying them.
\begin{teo}\label{teo_geo}
If $(X,d)$ is Polish and geodesic, then $(\P_2(X), W_2)$ is geodesic too. Furthermore,
consider $\mu_0,\mu_1\in \P_2(X)$ and a path
$t\mapsto \mu_t\in\P_2(X)$ from $\mu_0$ to $\mu_1$,
then $\mu_t$ is a constant speed geodesic
if and only if there exists $\mmu\in\P_2(\Geod(X))$ such that
$(e_0,e_1)_\#\mmu\in\Opt(\mu_0,\mu_1)$ and
$\mu_t=(e_t)_\#\mmu$.
\end{teo}

We consider on $X$ the Lebesgue measure $dk$. 
The measures $\mu$ that we are going to treat, are always induced by 
density  a.e.-continuous functions.

\begin{rmk}\label{rmk_et}
We know that the geodesic relying two points on a complete
riemannian manifold $X$ is almost everywhere unique (see for example \cite{villa08}).
This means that the set of pairs $(x,y)\in X^2$ such that the geodesic between them is unique, has
full measure. 
The same is true for any measure $\mu$ on $X$ if it is induced by a density a.e.-continuous function.

Therefore the maps
$e_t$ naturally induces almost everywhere a map that we denote in the same way:
$$e_t\colon X^2\to X,$$
sending the pair $(x,y)$ to the point $\gamma_t$ where $\gamma$ is
the geodesic such that $\gamma_0=x$ and $\gamma_1=y$.

Consider two probability measures $\mu_0,\mu_1$ over $X$ respecting the hypothesis of Theorem~\ref{teo_oneT}.
Let $T$ be the transport map between them,
then as  a consequence of Theorem \ref{teo_geo} the (unique, coming from the uniqueness of $T$) 
geodesic between $\mu_0$ and $\mu_1$ can be written as
$$\mu_t=(e_t\circ(\id,T))_\#\mu_0.$$
In what follows we will use the notation $e_t^{(T)}:=e_t\circ(\id,T)$,
and therefore we will have $\mu_t=(e_t^{(T)})_\#\mu_0$.\newline
\end{rmk}

\section{Reconstructing the visual input via the Gabor frame}\label{recon}

We are going to introduce the notion of continuous frame, this
allows the reconstruction of a function $I$ on the retinal 
plane, from the datum of an output function $\mu$.

\begin{defin}
Consider a Hilbert space $\H$ and a measure space $M$ with a positive measure~$\rho$.
A continuous frame is a family of vectors $\{\psi^k\}_{k\in M}$ such that $k\mapsto \langle f, \psi^k\rangle$
is a measurable function on $M$ for any $f\in \H$, and there exists $A,B>0$ such that
$$A\cdot \norm{f}^2\leq\int_X|\langle f,\psi^k\rangle|^2d\rho(k)\leq B\cdot \norm{f}^2.$$
\end{defin}

Consider $f,g\in \H$ and the mapping
$$h_f\colon g\mapsto \int_M\langle f,\psi^k\rangle\langle\psi^k,g\rangle d\rho(k).$$
This map is conjugated linear and moreover it is bounded. Indeed,
\begin{equation}\label{diseq}
|h_f(g)|^2\leq \int_M|\langle f,\psi^k\rangle|^2d\rho(k)\cdot\int_X|\langle \psi^k,g\rangle |^2d\rho(k)\leq B^2\norm{f}^2\norm{g}^2.
\end{equation}
By the Riesz' representation theorem, there exists a unique element $\underline h\in \H$ which verifies
 $h_f=\langle \underline h,-\rangle$. We denote this element by $\int_X\langle f,\psi^k\rangle \psi^k d\rho(k)$.

We denote by $S\colon \H\to \H$ the operator 
$$Sf:=\int_M\langle f,\psi^k\rangle\psi^kd\rho(k),\ \ \forall f\in \H.$$

\begin{lemma}\label{Sop}
the operator $S$ is linear and 
\begin{enumerate}
\item it is bounded and positive, with $\norm{S}\leq B$;
\item it is invertible;
\item the family $\{S^{-1}\psi^k\}_{k\in M}$ is a continuous frame;
\item for any $f\in \H$,
$$f=\int_M\langle f,\psi^k\rangle S^{-1}\psi^kd\rho(k)=\int_M\langle f,S^{-1}\psi^k\rangle \psi^kd\rho(k),$$
where the equality is intended in the weak sense.
\end{enumerate}
\end{lemma}
For a proof of this see \cite[\S5.8]{chris02}.\newline

From now on, we suppose $M=\R^2\times S^1\times \R^+$ with the measure $\rho$ 
such that  $d\rho(k)=\frac{dk}{\sigma^2}$.
In our case, the vectors $\psi^k$ are the Gabor filters
indexed by the points $k=(x,y,\theta,\sigma)\in M$. 
Consider the usual scalar product $\langle , \rangle$ in $L^2(\R^2)$, then
in Appendix \ref{appb} we prove that
there exists a constant $C_\psi\in \R^+$
such that for any pair of inputs $I,I'$ in $L^2(\R^2)$
\begin{equation}\label{eq_eq0}
\int_M\langle I,\psi^k\rangle\langle \psi^k,I'\rangle \frac{dk}{\sigma^2} 
=C_\psi\cdot\langle I,I'\rangle.
\end{equation}
As a corollary,
\begin{equation}\label{eq_eq}
SI=C_\psi\cdot I,
\end{equation}
 where the equality as to be intended in the weak sense, that is
$h_I=\langle C_\psi I,-\rangle_{L^2}$.
In Appendix~\ref{appc} we prove that the equality (\ref{eq_eq})
is also true in a much stronger sense.\newline

We observe that $\langle I,\psi_{x,y,\theta,\sigma}\rangle$ is exactly the output
function $\mu$ associated to $I$, see Definition~\ref{out}. Therefore by Lemma \ref{Sop} and equality (\ref{eq_eq}),
\begin{equation}\label{inv}
I(\tilde x,\tilde y)=\frac{1}{C_\psi}\cdot\int_M\mu(x,y,\theta,\sigma)\cdot\psi^k(\tilde x,\tilde y)
\cdot \frac{1}{\sigma^2}dxdyd\theta d\sigma.\newline
\end{equation}

In the next section, starting from two input
functions $I_0, I_1$ 
we produce a path $\mu_t$ of output functions.
In order to produce a path in the input space
from $\mu_t$,
we define
\begin{equation}\label{eq_imp0}
I_t:=\frac{1}{C_\psi}\cdot \int_M\mu_t(k)\cdot\psi^k \frac{dk}{\sigma^2}.
\end{equation}\newline

\section{Deformation of the output}\label{defoutput}
In this section we show the existence of a path relying
two output functions $\mu_0,\mu_1$. In particular for any output
we obtain two
probability measures from the positive and negative part of the output.
Using the results of Section \ref{optimal} and \ref{sec_geo} we 
build the paths relying the associated measure, and conclude from there.\newline

We introduce a new condition on the input image $I\colon \R^2\to [0,1]$.
If $\mu$ is the output associated to $I$, that is
$$\mu(k)= \langle I,\psi^k\rangle.$$
We look at $\mu$ as a function defined on the whole $X=\R^2\times S^1\times \R$
but it is supported only on $M=\R^2\times S^1\times \R^+\subset X$.
We define the two functions $\tilde\mu^+:=\max(\mu,0)$ and $\tilde\mu^-:=\max(-\mu,0)$.
We impose the condition
\begin{equation}\label{condition2}
\int_Xd^2(k, 0)\tilde \mu^+(k) dk< \infty,
\end{equation}
and the same for $\tilde \mu^-$.

\begin{lemma}\label{lemcond}
If conditions (\ref{condition2}) holds, then $\int d^2(k,k_0)\tilde\mu^+(k) dk$ is finite
for any $k_0\in X$.
\end{lemma}
\proof We start by observing that $\int_X\tilde \mu^+(k) dk<\infty$. Indeed,
if we consider the compact $B=\{k|\ d(k,0)\leq 1\}$ and the maximum $b$ of $\tilde \mu^+$ over $B$, 
we have
\begin{align*}
\int_X\tilde \mu^+(k)dk&\leq \tilde\mu^+(B)\cdot b+\int_{k\notin B} \tilde \mu^+(k)dk\\
&\leq \tilde \mu^+(B)\cdot b+\int_{k\notin B}d^2(k,0)\tilde\mu^+(k)dk
<\infty.
\end{align*}
For any $k_0\in X$, by the triangular inequality we have
$$d^2(k,k_0)\leq d^2(k,0)+d^2(k_0,0)+2\cdot |d(k,k_0)\cdot d(k_0,0)|\ \ \forall k\in X.$$
Therefore
\begin{align*}
\int_X d^2(k,k_0)\tilde \mu^+(k)dk&\leq \int_Xd^2(k,0)\tilde \mu^+(k)dk+d^2(k_0,0)\cdot\int_X\tilde \mu^+(k)dk+\\
&+2\cdot \int_X d(k,0)\cdot d(k_0,0)\tilde\mu^+(k)dk.
\end{align*}
The first two terms are clearly finite. Concerning the last one,
\begin{align*}
\int_X d(k,0)\cdot d(k_0,0)\tilde\mu^+(k)dk&\leq d(k_0,0)\cdot\tilde\mu^+(B)+\int_{k\notin B} d(k,0)\cdot d(k_0,0)\tilde \mu^+(k)dk\\
&\leq d(k_0,0)\cdot\tilde\mu^+(B)+d(k_0,0)\cdot \int_{k\notin B} d(k,0)^2\tilde \mu^+(k)dk<\infty,
\end{align*}
and this concludes the proof.\fine

\begin{rmk}\label{rmk_mu}
For $I\colon \R^2\to [0,1]$ we defined above
 $\tilde\mu^+$ and $\tilde\mu^-$.
Consider the coefficient
$$
m:=\int_M\tilde\mu^+(k)dk=-\int_M\tilde\mu^-(k)dk,
$$
which is well defined (see Remark \ref{rmk_mass})
and finite as a consequence of condition (\ref{condition2}) (see the proof of Lemma \ref{lemcond}).
We renormalize, $\mu^+:=\tilde\mu^+\slash m$ and $\mu^-:=\tilde\mu^-\slash m$.

Therefore for any function $\mu\colon M\to \R$
there exists two probability densities
$\mu^+$ and $\mu^-$ such that
$$\mu=m\cdot (\mu^+- \mu^-),$$
where $m$ is a positive coefficient and $\mu^+\cdot\mu^-\equiv 0$.
This means that $\mu^+$ and $\mu^-$ are the positive and negative part of $\mu$,
renormalized in order to become probability densities.
\end{rmk}

Given two inputs $I_0,I_1$, we define respectively $\mu_i^+,\mu_i^-$ and $m_i$ for $i=0,1$.
Using the equality~(\ref{inv}) we obtain,
\begin{align*}
I_0&=C_\psi^{-1}\cdot m_0 \cdot \left(\int_M\mu_0^+(k)\psi^k\frac{dk}{\sigma^2}-\int_M\mu_0^-(k)\psi^k\frac{dk}{\sigma^2}\right)\\
I_1&=C_\psi^{-1}\cdot m_1 \cdot \left(\int_M\mu_1^+(k)\psi^k\frac{dk}{\sigma^2}-\int_M\mu_1^-(k)\psi^k\frac{dk}{\sigma^2}\right).\\
\end{align*}

We consider the complete riemannian variety $X=\R^2\times S^1\times \R$ with
the Lebesgue metric~$dk$. In particular we consider $\mu_i^+$ and $\mu_i^-$, for $i=0,1$,
as measures on $X$ even if they are defined
over $M\subset X$. The function $\mu_i^\pm$ is identified
with $\mu_i^\pm\cdot dk$ if $\sigma>0$ and with the null measure elsewhere.
As a consequence of condition~(\ref{condition2}), $\mu_i^+,\mu_i^-\in \P_2(X)$ for $i=0,1$.

We consider the pairs $\mu_0^+,\mu_1^+$ and $\mu_0^-,\mu_1^-$,
by Theorem~\ref{teo_geo} and
 Remark \ref{rmk_et}, there exists two constant speed geodesics
$\mu_t^+$ and $\mu_t^-$.

\begin{rmk}\label{mut}
In particular by Theorem \ref{teo_oneT} there exists a transport map $T^+$ such that $\mu_1^+=T^+_\#\mu_0^+$ and
the same for the negative part with a transport map~$T^-$. Therefore,
$\mu_t^+=\left(e_t^{(T^+)}\right)_\#\mu_0^+$ and $\mu_t^-=\left(e_t^{(T^-)}\right)_\#\mu_0^-$.
\end{rmk}

\begin{rmk}\label{mut2}
By definition of $e_t$ and of the transport map, the measures $\mu_t^\pm$
are null outside $M$ for any $t\in [0,1]$, therefore we can always
look at them as measures in $\P_2(M)$.

Furthermore, we point out that by construction $\mu_t^+$ and $\mu_t^-$ are 
absolutely continuous with respect to the Lebesgue measure $dk$.
In the following
we will use the notation $\mu_t^+$ and $\mu_t^-$
indistinctly for the measures and for the density functions (defined over $M$)
when there is no risk of confusion.\newline
\end{rmk}

We consider a linear variation of the mass $m$,
this means that we define a varying coefficient
\begin{equation}\label{varco}
m_t:=m_0(1-t)+m_1t\ \ \forall t\in [0,1]
\end{equation}
We define the path of output functions using these coefficients,
\begin{equation}\label{path}
\mu_t:=m_t\cdot\left(\left(e_t^{(T^+)}\right)_\#\mu_0^+-\left(e_t^{(T^-)}\right)_\#\mu_0^-\right).
\end{equation}
Finally from Equation (\ref{eq_imp0}) we 
obtain a path of input functions from $I_0$ to $I_1$.\newline

\section{Constraining the output}\label{const}

In this section we introduce a useful tool to describe the 
geodesic $\mu_t$,
the so called weak riemannian structure of $(\P_2(X), W_2)$,
the space of probabilities endowed with the Wassertein distance.

If $\mu_t$ is an absolutely continuous curve in $\P_2(X)$
(with respect to the Wassertein distance), consider a time dependent
vector field $v_t$ on $TX$ such that the following \cor{continuity equation}
is verified in the sense of distributions,
\begin{equation}\label{eq_cont}
\frac{d}{dt}\mu_t+\nabla\cdot(v_t\mu_t)=0.
\end{equation} 
For the proof of the following theorem we refer again to \cite{amgi13}.
\begin{teo}[see {\cite[Theorem 2.29]{amgi13}}]\label{teoab}
If $X$ is a smooth complete riemannian manifold without boundary, then
\begin{enumerate}
\item for every absolutely continuous curve $\mu_t\in \P_2(X)$ there exists a Borel family
of vector fields $v_t$ such that $\norm{v_t}_{L^2(\mu_t)}\leq |\dot \mu_t|$ for a.e.~$t$ and
the continuity equation (\ref{eq_cont}) is satisfied (in the sense of distributions);
\item if $(\mu_t,v_t)$ satisfies (\ref{eq_cont}) and $\int_0^1\norm{v_t}_{L^2(\mu_t)}dt$ is finite,
then $\mu_t$ is an absolutely continuous curve (up to a negligible set of points)
and $|\dot\mu_t|\leq \norm{v_t}_{L^2(\mu_t)}$ for a.e.~$t\in[0,1]$.
\end{enumerate}
\end{teo}

\begin{rmk}
We recall the Benamou-Brenier formula proved for example
at \cite[Proposition 2.30]{amgi13} and stating that the minimization problem
solved by a geodesic relying $\mu_0,\mu_1\in\P_2(X)$ can be reformulated
in terms of the vector field $v_t$. In particular we have
$$W_2(\mu_0,\mu_1)=\inf \int_0^1\norm{v_t}_{L^2(\mu_t)}dt,$$
where the infimum is taken among all weakly continuous distibutional solutions of the continuity equation 
for $(\mu_t,v_t)$.
\end{rmk}

As a direct consequence of Theorem \ref{teoab}, for every absolutely continuous curve $\mu_t$ in $\P_2(X)$,
there exists a family of vector fields $(v_t)$ verifying the continuity equation
and such that $\norm{v_t}_{L^2(\mu_t)}=|\dot\mu_t|$ for a.e.~$t$.
This family is not unique in general, but it is unique if we define
as follows the tangent space to $\P_2(X)$ where the vector fields
must live in.
\begin{defin}
If $\mu\in\P_2(X)$ then the tangent space to $\P_2(X)$ is defined as
\begin{align*}
T_\mu\P_2(X)&:={\overline{\left\{\nabla \varphi:\ \varphi\in \cc^\infty_c(X)\right\}}}^{L^2(\mu)}\\
&=\left\{v\in L^2(\mu):\ \int \langle v,w\rangle d\mu=0,\ \forall w\in L^2(\mu)\ \m{s.t.}\ \nabla\cdot(w\mu)=0\right\}.
\end{align*}
\end{defin}

Therefore for any absolutely continuous curve $\mu_t$ in $\P_2(X)$,
we have an associated vector field $v_t$. In particular we have
it for the geodesics $\mu_t^+$ and $\mu_t^-$
obtained via two inputs $I_0,I_1$ (see Remark \ref{mut}).
Observe that both $\mu_t^+$ and $\mu_t^-$ are constant
speed geodesics in $(\P_2(X),W_2)$ therefore
they are absolutely continuous curves, and so
Theorem \ref{teoab} applies to them.

We denote by $v_t^+$ and $v_t^-$ the vector fields associated respectively
to $\mu_t^+$ and $\mu_t^-$. Moreover, we define the normalized image
$$J_t:=\frac{I_t}{m_t},$$
so that
\begin{equation}\label{eq_0}
J_t=\int_M\frac{\mu_t(k)}{m_t}\psi^k\frac{dk}{\sigma^2}=\int_M(\mu^+_t-\mu^-_t)\psi^k\frac{dk}{\sigma^2}.
\end{equation}

We are interested in the existence of a family $J_t$ of inputs $\R^2\to [0,1]$
that relies $J_0$ to $J_1$, such that $J_t$ is in $\cc^1([0,1]; L^2(\R^2))$
and $\frac{\mu_t(k)}{m_t}=\langle J_t,\psi^k\rangle $ for any $k\in M$.

In order to find such a path, we observe that if it exists, then
$$\frac{d}{dt}\frac{\mu_t(k)}{m_t}=\left\langle \frac{dJ_t}{dt},\psi^k\right\rangle \ \ \forall k\in M\subset X.$$
From the continuity equation we know that
\begin{align*}
\frac{d\mu_t^+}{dt}&=-\nabla\cdot(v_t^+\mu^+_t)\\
\frac{d\mu_t^-}{dt}&=-\nabla\cdot(v_t^-\mu^-_t).
\end{align*}
We define $v_t:=v_t^+-v_t^-$, and therefore
$$\frac{d}{dt}\frac{\mu_t}{m_t}=-\nabla\cdot\left(v_t\frac{\mu_t}{m_t}\right).$$
Indeed, by Theorem \ref{teoab}, $v_t^+\in L^2(\mu_t^+)$ and $v_t^-\in L^2(\mu_t^-)$.
Therefore it is possible to extend both vector fields to the whole $M$ in such a way that $v_t^+\cdot\mu_t^-\equiv 0$
and $v_t^-\cdot \mu_t^+\equiv 0$.
Then we have,
$$\frac{d}{dt}\frac{\mu_t(k)}{m_t}=-\nabla\cdot \left(v_t\langle J_t,\psi^k\rangle\right)=-\left\langle J_t,v_t(k)\cdot\nabla\psi^k+\psi^k\cdot(\nabla\cdot v_t(k))\right\rangle
\ \ \ \forall k\in M\subset X.$$

For an opportune vector field $\alpha\in TM$ we have 
\begin{equation}\label{alpha}
\nabla\psi^k=\psi^k\cdot \alpha(k).
\end{equation}
In particular if we use the notation $(\tilde x_k,\tilde y_k)=\sigma^{-1}R_{-\theta}(\tilde x-x,\tilde y-y)$,
where $k=(x,y,\theta,\sigma)\in M$,
it is straightforward to verify that
\begin{align*}
\alpha^x&=2 \left(\sigma^{-2}(\tilde x-x)+\frac{\sigma^{-1}\sin\theta}{\tan(2\tilde y_k)}\right)\\
\alpha^y &= 2\left(\sigma^{-2}(\tilde y-y)-\frac{\sigma^{-1}\cos\theta}{\tan(2\tilde y_k)}\right)\\
\alpha^\theta &= \frac{2\sigma^{-1}}{\tan(2\tilde y_k)}\left(\cos\theta(\tilde x-x)+\sin\theta(\tilde y-y)\right)\\
\alpha^\sigma &= 2 \left(\sigma^{-3}\norm{\tilde v-v}^2+\frac{\sigma^{-2}(\sin\theta(\tilde x-x)-\cos\theta(\tilde y-y))}{\tan(2\tilde y_k)}- \frac{3\sigma^{-1}}{4}\right).
\end{align*}
In particular for any $k$, the vector field $\alpha$ is
well defined almost everywhere on $\R^2$.

Therefore for a.e.~$t\in [0,1]$ and every $k\in M$,
$$\left\langle \frac{dJ_t}{dt},\psi^k\right\rangle=-\left\langle J_t(\alpha\cdot v_t+\nabla\cdot v_t),\psi^k\right\rangle.$$
Thus, if $\alpha\cdot v_t+\nabla\cdot v_t$ is independent  
of the variable $k$ as a function $\R^2\to \R$, 
as $\left(\psi^k\right)_k$ is a frame, this implies 
\begin{equation}\label{eq_a}
\frac{dJ_t}{dt}=-J_t(\alpha\cdot v_t+\nabla\cdot v_t).
\end{equation}
This last equality is true in the weak sense but also
in the (stronger) sense showed in Appendix~\ref{appc}.\newline

In order to state our last theorem,
we introduce two additional conditions. First
we impose that the inputs $I_0,I_1$ are null outside a compact
subset of $\R^2$. This simply says
that the images we are treating are limited in space.

We also impose that the vector fields take their values
in a Sobolev space.

\begin{teo}
Consider two inputs $I_0,I_1\colon \R^2\to [0,1]$ in $L^2(\R^2)$ null outside a compact
subset of $\R^2$, their associated output functions $\mu_0,\mu_1$, 
 the absolutely continuous curve $\mu_t$ relying $\mu_0$ to $\mu_1$ defined in~(\ref{path})
 and the associated (unique) vector field $v_t$.
 Moreover, take the vector field $\alpha$ defined in (\ref{alpha}).
 
 We suppose that 
 $$v\in L^1([0,1]; W^{1,\infty}(X,\R^4)).$$
If for any $t\in [0,1]$ the following equality is verified
\begin{equation}\label{eq_b}
\nabla\left(\alpha\cdot v_t+\nabla\cdot v_t\right)=0,
\end{equation}
then there exists a path $I_t\in\cc^1([0,1];L^2(\R^2))$
relying $I_0$ to $I_1$
such that $\mu_t(k)=\langle I_t,\psi^k\rangle$ for any $k\in M$.
\end{teo}
\proof We use the notation $u_t:=\alpha\cdot v_t+\nabla\cdot v_t$
and observe that by (\ref{eq_b}), $u_t$ is independent of the point $k$.
As a consequence of this independence, $u_t$ can
be defined also where $\al$ is singular.
We consider the differential equation
$$\frac{dJ_t}{dt}=-u_t\cdot J_t,$$
and the coefficient $m_t$ as defined in (\ref{varco}).
Given the initial function $J_0:=I_0\slash m_0$, the solution to the equation above is
$$J_t:=J_0\cdot e^{-h_t},$$
where $h_t:=\int_0^tu_sds\ \forall t\in[0,1]$ is a primitive of $u_t$,
and it exists as a consequence of $v_t\in W^{1,\infty}$.
The function $J_0$ is null outside a compact subset $K\subset \R^2$,
and $h_t$ is continuous therefore limited over $K$. This implies $J_t\in L^2(\R^2)$
for any $t\in [0,1]$.\newline

We define $I_t:=m_t\cdot J_t$ and
$$\nu_t(k):=m_t\cdot \langle J_t,\psi^k\rangle=\langle I_t,\psi^k\rangle.$$
By construction, for any $t\in [0,1]$ $\nu_t$ satisfies a.e.~the continuity equation (\ref{eq_cont})
with respect to the vector field $v_t$. Moreover,
 $\nu_0=\frac{\mu_0}{m_0}=\mu^+_0-\mu^-_0$.
The solution to the continuity equation under this conditions
is unique for absolutely continuous measures (see \cite{amcri08}).
Therefore $\nu_t=\frac{\mu_t}{m_t}$ a.e.~for a.e.~$t\in[0,1]$.
We observe that this also proves that $m_1J_1=I_1$.\fine

\section{Discrete model}\label{discmod}

In order to produce an implementation of a transport model
in the space $M=\R^2\times S^1\times \R^+$,
we have to work in a discrete setting, and therefore
it is not possible to use the frame $\psi^k$ introduced above,
because it is in fact a frame only when the index $k$ varies
along all $M$.

We consider instead a discrete frame that produce a new set
of output functions. The frame we consider is the so called Gabor Wavelet Pyramid
(see for instance \cite{knpg08}) that relies again in the Gabor 
complex mother function
$$e^{-\tilde x^2-\tilde y^2}\cdot e^{2\pi i\omega \tilde y},$$
but we split it in its real and complex component, defining two
sets of filter functions associated
to the two following mother functions
\begin{align*}
\psi_e(\tilde x,\tilde y)&:=e^{-\tilde x^2-\gamma\cdot \tilde y^2}\cdot \cos(2\pi \omega\tilde y)\\
\psi_o(\tilde x,\tilde y)&:=e^{-\tilde x^2-\gamma\cdot \tilde y^2}\cdot \sin(2\pi\omega\tilde y).
\end{align*}

As we have seen, in the case of the continuous frame, the elements of the frame
are obtained via the action of the Heisenberg group on the 
mother function. In the case of wavelets the acting group is the affine group.
We recall that by $\psi_{e,\theta}$ and $\psi_{o,\theta}$ we mean the 
same functions with the rotation by $\theta$ applied, that is
\begin{align*}
\psi_{e,\theta}(\tilde x,\tilde y)&=\psi_e\left(R_{-\theta}(\tilde x,\tilde y)\right)\\
&=\psi_e\left(\cos(\theta)\tilde x+\sin(\theta)\tilde y,\ -\sin(\theta)\tilde x+\cos(\theta)\tilde y\right),
\end{align*}
and analogously for $\psi_o$.

\begin{defin}\label{def_wavelets}
Given two real numbers $a_0,b_0$ and a positive integer number $d$,
we consider $\theta_0:=\pi\slash d$ and $\theta_\ell=\ell\cdot \theta_0$ for any $\ell$
positive integer. Then we can define the wavelets associated to $\psi_e$ and $\psi_o$.
For any $n,k,\ell,j\in \Z_{\geq0}$,
\begin{align*}
\psi_e^{n,k,\ell,j}(\tilde x,\tilde y)&=\frac{1}{a_0^j}\cdot \psi_{e,\theta_\ell}\left(\frac{\tilde x}{a_0^j}-nb_0,\ \frac{\tilde y}{a_0^j}-kb_0\right)\\
\psi_o^{n,k,\ell,j}(\tilde x,\tilde y)&=\frac{1}{a_0^j}\cdot \psi_{o,\theta_\ell}\left(\frac{\tilde x}{a_0^j}-nb_0,\ \frac{\tilde y}{a_0^j}-kb_0\right)
\end{align*}

\end{defin}

As developed in \cite{daub90} by Daubechies for one variable wavelets, and generalized
to two variable wavelets by Lee in \cite{lee96}, for an opportune choice of $a_0,b_0\in \R$,
and $n,k,\ell,j\in \Z_{\geq0}$, the wavelets above form a frame, that is there exist real numbers 
$A,B>0$ such that for any $f\in L^2(\R)$ we have
$$A\cdot \norm{f}^2\leq \sum_{n,k,\ell,j}| \langle f,\psi_e^{n,k,\ell,j}\rangle|^2+\sum_{n,k,\ell,j}|\langle f,\psi_o^{n,k,\ell,j}\rangle|^2
\leq B\cdot \norm{f}^2,$$
where the scalar product is the classical
$$\langle f,\psi\rangle=\int f(\tilde x,\tilde y)\cdot \psi(\tilde x,\tilde y)d\tilde xd\tilde y.$$
As a consequence it is possible to reconstruct the function, as
\begin{equation}
f=C\cdot \left(\sum_{n,k,\ell,j}\langle f,\psi_e^{n,k,\ell,j}\rangle\cdot \psi_e^{n,k,\ell,j}+\sum_{n,k,\ell,j}\langle f,\psi_o^{n,k,\ell,j}\rangle
\cdot\psi_o^{n,k,\ell,j}\right),
\end{equation}
where $C$ is a constant and the equality is true in the weak sense.
Observe that the splitting of the mother function in its even and odd part allows
to consider only `half' of the $S^1$ circle, or equivalently
it allows to use an equipartition of the space of directions $S^1\slash\{\pm1\}$.\newline

In our setting, we consider in this case four output functions obtained
by convoluting any input $I\colon \R^2\to [0,1]$ with the 
frame above. Then, we apply the same procedure of Section \ref{defoutput},
that is we use optimal transport between functions on $M$ to
build a path relying the initial output to the final output.\newline

Instead of working over the whole $M$, we consider a compact subset
$$M_c:=[0,D]^2\times S^1\times [\sigma_{\min},\sigma_{\max}],$$
where $[0,D]^2$ represents the portion of the plane which is 
registered on the retina, while $[\sigma_{\min}, \sigma_{\max}]$
is the interval where the scale parameter $\sigma$ varies.
 
 \begin{rmk}\label{gwp}
The idea behind the Gabor Wavelet Pyramid is that we define the 
output functions indexed by a discrete subset of $M_c$. We use the notation  $\llbracket0,N\rrbracket$ for
the subset $\{0,1,\dots,N\}\subset \Z_{\geq0}$ for any positive integer $N$. We set $a_0=\sigma_{\min}$ and
$j_\star=\log_{a_0}(\sigma_{\max})$. For any $j\in \llbracket1,j_\star\rrbracket$
  we consider the discrete subset
$$(b_0a_0^j)\cdot \left\llbracket0,\left\lfloor \frac{D}{b_0a_0^j}\right\rfloor\right\rrbracket^2\times \theta_0\cdot [1,d]\times \left\{a_0^j\right\}\subset
[0,D]^2\times S^1\times \left\{a_0^j\right\},$$
thus defining different strata of a discrete subset of the whole $M_c$.
\end{rmk}

\begin{rmk}\label{rmkwav}
We observe that
any wavelet $\psi_e^{n,k,\ell,j}$ (the same is true for $\psi_o$)
is centered at $(nb_0a_0^j,\ kb_0a_0^j)$. Furthermore,
given two inputs $I_0,I_1$, 
we can define the following functions on $M_c$, for
$i=0,1$,
\begin{align*}
\tilde \mu_{i,e}^+(nb_0a_0^j,kb_0a_0^j,\theta_\ell,a_0^j)&=\max\left(\langle I_i,\psi_e^{n,k,\ell,j}\rangle,0\right)\\
\tilde \mu_{i,e}^-(nb_0a_0^j,kb_0a_0^j,\theta_\ell,a_0^j)&=-\min\left(\langle I_i,\psi_e^{n,k,\ell,j}\rangle,0\right)
\end{align*}
and analogously for $\tilde \mu_{i,o}^+$ and $\tilde \mu_{i,o}^-$.
We normalize these functions in order to obtain the probabilities $\mu_{i,e}^{\pm}$
 for $i=0,1$ (and analogously for $\mu_{i,o}^{\pm}$).\newline
 \end{rmk}

In order to build the transport path $\mu_t$ for any of the pairs of output functions above,
we focus on the distance over $M_c$. In fact, 
instead of working with the distance treated in the previous sections,
another distance turns out to be more efficient in our setting,
and at the same time it is a distance equivalent to the previous one.

\begin{defin}
Two distances $d_1,d_2$ over the space $\Omega$ are said to be equivalent
if for any compact set $K\subset\Omega$, there exists a constant $C$ such that
$$\frac{1}{C}\cdot d_2(p_1,p_2)\leq d_1(p_1,p_2)\leq C\cdot d_2(p_1,p_2)\ \ \forall p_1,p_2\in \Omega.$$

They are locally equivalent if for any $p_0\in\Omega$ there exists a neighborhood $U$
of $p_0$ such that $d_1$ and $d_2$ are equivalent over $U$.
\end{defin}

The main reference for the equivalence result we are going
to use is \cite{nsw85}.
With our notation (see Remark \ref{rmk_met}),
we consider in every point the metric
$$\tilde g(x,y,\theta,\sigma):=
\left(\begin{array}{cccc}
1 &0 &0 &0\\
0 &1 &0&0\\
0 &0 &\frac{h_1^2}{\sigma^2} &0\\
0& 0& 0& \frac{h_2^2}{\sigma^2}\\
\end{array}\right)
$$
and we denote by $Y_1,Y_2,Y_3,Y_4$ the orthonormal basis at every point
defined by $Y_1=X_1$ and $Y_2=X_2$, while $Y_3=\frac{\sigma X_3}{h_1}$
and $Y_4=\frac{\sigma X_4}{h_2}$.
As we are interested in constant coefficient flows, consider
$$Y=c_1Y_1+c_2Y_2+c_3Y_3+c_4Y_4,$$
with $c_i\in \R^+$ for $i=1,\dots,4$.
In order to evaluate the distance $d_c$ between two points $p_0=(x_0,y_0,\theta_0,\sigma_0)$
and $p_1=(x_1,y_1,\theta_1,\sigma_1)$, we consider the constant coefficient vector field $Y$
that induces a curve $p_t\colon [0,1]\to M$ relying them. We develop these evaluations
 in Appendix~\ref{appd}.

\begin{prop}[See, {\cite[Theorem 2]{nsw85}}]
The distance induced by the riemannian metric $\tilde g$ and the distance $d_c$
are locally equivalent.
\end{prop}

In the next section we are going to show an implementation
of our transport model based on the distance $d_c$.

\begin{rmk}\label{rmkit}
From the previous results we know that there exists four unique
paths $\mu_{t,e}^{\pm},\ \mu_{t,o}^{\pm}$ relying the associated probabilities.
Using the fact that the chosen set of filters is a frame, we reconstruct
the intermediary images, by defining
\begin{align*}
I_t:=&\sum_{n,k,\ell,j}m_{e,t}(\mu_{t,e}^+-\mu_{t,e}^-)(nb_0a_0,kb_0a_0,\theta_\ell,a_0^j)\cdot\psi_e^{n,k,\ell,j}\\
&+\sum_{n,k,\ell,j}m_{o,t}(\mu_{t,o}^+-\mu_{t,o}^-)(nb_0a_0,kb_0a_0,\theta_\ell,a_0^j)\cdot\psi_{o}^{n, k,\ell,m}.
\end{align*}
Where the coefficients $m_{e,t}$ and $m_{o,t}$ depend
on the masses of $\tilde \mu_{i,e}^\pm$ and $\tilde \mu_{i,o}^\pm$
for $i=0,1$.\newline
\end{rmk}

\section{Implementation}\label{impl}
In order to implement our model, we code a Gabor Wavelet Pyramid as sketched in
the previous section, and evaluate the transport maps
using a Sinkhorn's algorithm of the kind treated by Peyré and Cuturi \cite{peycut19}.
In discrete setting, the transport plan between two functions $u_1,u_2$,
is a matrix $P$. In particular, if the two functions are represented as vectors
of dimension~$m$, the matrix $P$ has dimension $m\times m$ and verifies the property
\begin{equation}\label{cond}
P\bm{1}_m=u_1,\ \ P^T\bm{1}_m=u_2
\end{equation}
where $\bm{1}_m$ is the vector of length $m$ whose coordinates are all $1$,
and therefore (\ref{cond}) means that $P\in \adm(u_1,u_2)$.

We recall that the transport plan $P$ between two functions $u_1,u_2$
minimizes the product $\langle P,C\rangle$ where $C$ is the cost matrix associated
to our setting, and $P$ respects condition (\ref{cond}).
In the case of the Sinkhorn's algorithm, with regularization coefficient $\varepsilon$,
the minimized quantity is
\begin{equation}\label{condeq}
\langle P,C \rangle -\varepsilon H(P),
\end{equation}
where $H(P)$ is the entropy of $P$. The uniqueness of the solution
is in fact true for any strictly concave function $H$ (see \cite{dimager20}
for a wider analysis on entropy regularizations).
A well known result \cite[\S4]{peycut19} states that if
$L_C(u_1,u_2)$ is the minimum reached by $\langle P,C\rangle$
under condition (\ref{cond}), and $L^\varepsilon_C(u_1,u_2)$
is the minimum of (\ref{condeq}) under the same condition,
then
$$L^\varepsilon_C(u_1,u_2)\xrightarrow{\varepsilon\to 0} L_C(u_1,u_2).$$
We observe that the plan achieving $L_C(u_1,u_2)$ is not always unique
in the discrete setting, while the one achieving $L_C^\varepsilon(u_1,u_2)$
is in fact unique as proved for example by \cite[Proposition 4.3]{peycut19}.
Moreover, if $P^\varepsilon$ is the optimal plan for $L^\varepsilon_C$, then
$P^\varepsilon\to P^\star$ where $P^\star$ is the maximum entropy plan
among those achieving the optimum $L_C$.\newline

In order to compute the geodesic between $\mu_{0,e}^+$ and $\mu_{1,e}^+$
we define the optimal plan for this case
$$P^+_e:=\argmin_{P\in \adm(\mu_{0,e}^+,\mu_{1,e}^+)}\left(\langle P,C\rangle -\varepsilon \cdot H(P)\right).$$
and consider the meshgrid of Remark \ref{gwp}.
The pseudocode for the computation is the following.

\begin{algorithm}[H]
        \begin{algorithmic}
	\caption{Computation of the geodesic $\mu_{t,e}^+$} \label{alg1}
          \STATE \textbf{Input}  the optimal plan $P_e^+$; $\sigma_{\min},\sigma_{\max},b_0\in \R^+$; $D,d\in \Z_{>0}$;
          \STATE Set $\mu_{t,e}^+\equiv 0$;
          \STATE Set $a_0=\sigma_{\min}$;
          \STATE Set $j_\star=\log_{a_0}(\sigma_{\max})$;
          \FOR{ $j,j'=1,\dots, j_\star$}
          \FOR{ $\ell,\ell'=1,\dots, d$}
          \FOR{ $n,k=0,\dots, \left\lfloor\frac{D}{b_0a_0^j}\right\rfloor$ and $n',k'=0,\dots,\left\lfloor\frac{D}{b_0a_0^{j'}}\right\rfloor$}
          \STATE Set $p,q=(nb_0a_0^j,kb_0a_0^j, \ell\theta_0, a_0^j),(n'b_0a_0^{j'},k'b_0a_0^{j'},\ell'\theta_0,a_0^{j'})$;
           \STATE Take the geodesic $\gamma$ relying $p$ to $q$ and its point $\gamma_t$;
           \STATE $\mu_{t,e}^+(\gamma_t):=\mu_{t,e}^+(\gamma_t)+P_e^+(p,q)$;
           \ENDFOR
           \ENDFOR
           \ENDFOR
           \STATE \textbf{Output} $\mu_{t,e}^+$.
        \end{algorithmic}
      \end{algorithm}    	

We do the analogously for $\mu_{t,e}^-$ and $\mu_{t,o}^\pm$.
In order to develop the code we implemented various algorithms developed
in \cite{peycut19}, in particular we used part of the code developed by the same
authors for the Sinkhorn's algorithm and disposable online \cite{peycut19.2}.\newline

In our simulation, we set $a_0=2, b_0=1$ and as already said $D=32$ or $64$.
We also set $\ell=8$ which is coherent with previous results on the 
orientation sampling in primates (see \cite{bcs14, ring02}).
Furthermore, we set $\sigma_{\min}=1,1244$ and $\sigma_{\max}=\sigma_{\min}\cdot D$. 
The pseudocode for the complete simulation is the following.

\begin{algorithm}[H]
        \begin{algorithmic}
	\caption{Computation of $I_t$} 
          \STATE \textbf{Input} $h_1,h_2, b_0,\sigma_{\min}, \sigma_{\max},\varepsilon\in \R^+$; $N_{\iter},D,d\in \Z_{>0}$; $I_0,I_1\colon \llbracket 1,D\rrbracket\to [0,1]$;
          \STATE Read the two input image $I_0,I_1$ that are $D\times D$ matrices (with usually $D=32,64$);
          \STATE Evaluation of the cost matrix $c$ (see Appendix \ref{appd}) over the discrete subset of $M_c$ built in Remark \ref{gwp};
          \STATE Evaluation of $\mu_{i,e}^{\pm}$ and $\mu_{i,o}^{\pm}$ (see Remark \ref{rmkwav});
          \STATE Evaluation of $P_{e}^\pm$ and $P_o^\pm$ via the Sinkhorn's algorithm with $N_{\iter}$ iterations and $\varepsilon$ regularization coefficient;
          \STATE Evaluation of $\mu_{t,e}^\pm$ and $\mu_{t,o}^\pm$ via Algorithm \ref{alg1};
          \STATE Evaluation of $I_t$ via the formula of Remark \ref{rmkit};
           \STATE \textbf{Output} $I_t$.
        \end{algorithmic}
      \end{algorithm}


The images we consider are simple shapes of the letters `T' and `E', 
the second ones rotated of about $\frac{\pi}{4}$ counterclockwise.
We show the implementation for this input in order to emphasize the ability of 
our model to reconstruct rotational displacement. We also add an hammer-type shape,
in order to show that the implementation works well also in the case
of a multi-scale object, that is shapes with different thickness along
the figure are also well preserved.

Furthermore, we
compare this numerical cortical-style implementation with a classical 2-dimensional planar
regularized optimal transport implementation, following 
the general theory treated for example in \cite[\S3]{peycut19},
and evaluating the transport matrix again via the Sinkhorn's algorithm.
In this case we take the two inputs $I_0,I_1$ and the probabilities
$\nu_0,\nu_1\in \P(\llbracket 1,D\rrbracket)$ obtained by normalizing the inputs. 
The cost function is the quadratic cost $c\colon \llbracket 1,D\rrbracket^2\to \R^+$ such that $c(i,j)=(i-j)^2$.\newline

We point out that we applied a smoothing threshold
in order to reduce the blur due to the entropic regularization
in both our implementations. In particular, we passed the pictures
through the sigmoid 
$$\Sigma(z)=\frac{1}{1+e^{-k\cdot (z-z_0)}}$$
where $y\in [0,1]$ is the pixel intensity and we choose
$k=30$ and $z_0=0.65$.\newline 



\begin{figure}[h]
\includegraphics{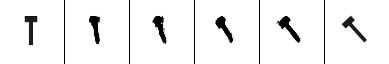}
\caption{Reconstruction of missing images with the parameters $D=64$, $h_1=0.7$, $h_2=5$,
 $N_{\iter}=1000$ and $\varepsilon=0.04$. We also processed
 the image applying a smooth threshold via the $\Sigma$ sigmoid function.}
\end{figure}

 \begin{figure}[h]
\includegraphics{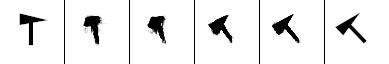}
\caption{Reconstruction of missing images in the case of the hammer shape, with the parameters $D=64$, $h_1=0.7$, $h_2=5$,
 $N_{\iter}=1000$ and $\varepsilon=0.04$. We also processed
 the image applying a smooth threshold via the $\Sigma$ sigmoid function.}
\end{figure}

\begin{figure}[h]
\includegraphics{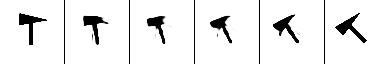}
\caption{Reconstruction of missing images via the application
of a planar optimal transport in the case of the hammer shape, with cost induced 
by the square of the euclidean distance, and using the
the Sinkhorn's algorithm with $N_{\iter}=10000$ and $\varepsilon=0.01$.
We also processed the image through the sigmoid $\Sigma$.}
\end{figure}

\begin{figure}[h]
\includegraphics{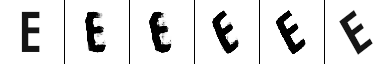}
\caption{Reconstruction of missing images with the parameters $D=64$, $h_1=0.7$, $h_2=5$,
 $N_{\iter}=1000$ and $\varepsilon=0.04$.
 We also processed the image through the sigmoid $\Sigma$.}
\end{figure}

\begin{figure}[h]
\includegraphics{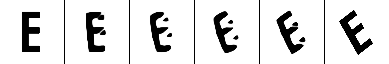}
\caption{Reconstruction of missing images via the application
of a planar optimal transport, with cost induced 
by the square of the euclidean distance, and using the
the Sinkhorn's algorithm with $N_{\iter}=10000$ and $\varepsilon=0.01$.
We also processed the image through the sigmoid $\Sigma$.}
\end{figure}

In these reconstructed images obtained with the cortical model,
we observe the conservation of the 
image structure along the rotational movement, while in the standard
2-d optimal transport the basic shape is lost. Therefore, lifting
the input via Gabor filters and moving the output function 
through optimal transport tools, seems to allow
an effective image deformation preserving the fundamental
aspects of rigid rotational motion.

\begin{figure}[h]
\includegraphics{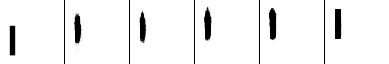}
\caption{Reconstruction of missing images in the case of the `I' shape,
 with the parameters $D=64$, $h_1=0.7$, $h_2=5$,
 $N_{\iter}=1000$ and $\varepsilon=0.01$.
 We also processed the image through the sigmoid $\Sigma$.}\label{ishp}
\end{figure}

If we consider a rigid translation, our procedure works
well in the case of simple shapes that
are moved along their principal direction, as it is
the case for the `I' shape in Figure \ref{ishp}.\newline

\section{Conclusions and future developments}

We formulated the general theory of a lifting of retinal inputs
in a 4-dimensional cortical space, and of the time completion between
two cortical outputs $\mu_0,\mu_1$ (corresponding
to inputs $I_0,I_1$ in the retinal plane). We did the lifting through Gabor filters, and we 
considered the frame conditions that allow to project cortical measures
back to the retinal space, thus obtaining from the completion paths
$\mu_t^{\pm}$ on the cortical space, a retinal path $I_t$ 
between the original inputs.

We obtained the cortical paths via methods of optimal transport,
where the cost function is the squared distance on the cortical space.
We implemented these tools using a Sinkhorn's algorithm
on a discretized version of the cortical manifold, and via a Gabor Wavelet Pyramid
system we also implemented the retinal path $I_t$. We tested positively
our model on rigid rotational movements of multi-scale shapes, verifying
the shape conservation.

The Sinkhorn's algorithm does not minimize the original cost $\langle c,\gamma\rangle$,
but a regularized version of it with an entropic correction $\langle c,\gamma\rangle -\varepsilon\cdot h(\gamma)$.
If $\gamma^\varepsilon$ is a solution to the regularized problem, we know that $\gamma^\varepsilon\xrightarrow{\varepsilon \to 0}\gamma^\star$
where $\gamma^\star$ is a solution of the original problem. 
The difference between $\gamma^\varepsilon$ and $\gamma^\star$ is what induces
a blurring on the resulting moving shapes in our implementation. We de-blurred the images by filtering
them through a sigmoid function. 

Knowing that $\norm{\gamma^\varepsilon-\gamma^\star}$
is controlled by $\varepsilon$, in future works
we intend to follow the Rigollet-Weed model \cite{rigoweed18} and use another
optimal transport technique to do the de-blurring. Indeed, we can suppose
that the noise factor $\sigma^2$ in \cite{rigoweed18} depends directly on $\varepsilon$,
and from there we can develop the Wasserstein distance minimization
as conceived in the work above.

In the future we also intend to improve the representation of translational movements. In our actual implementations,
 translation works correctly particularly along the boundary directions.
In order to extend our model by considering spatio-temporal Gabor filtering, 
as done for example in \cite{bccs14}, we could improve translation in the direction ortogonal to boundaries,  
still implementing the optimal transport tools considered in this work and thus preserving
 the boundary shapes of the retinal inputs. 
 
 Anyway, the problem with displacements orthogonal to the object boundaries 
 is linked also to the structure of the Gabor Wavelet Pyramid: in order to optimize the sampling, the pyramid frame is usually built by considering 
 $\sigma$ values that are ``far” from $0$, this emphasizes the boundary constrains. 
 We will work in frame buildings that allow for different $\sigma$ ranges.


\newpage
\newpage

\appendix
\section{}\label{appa}
In this first appendix, we prove that 
the integral 
\begin{equation}\label{int}
\int_{M}\psi^k(\tilde x,\tilde y)dk,
\end{equation}
is well defined for any pair $(\tilde x,\tilde y)\in \R^2$ and
has finite value $0$,
independently from the pair $(\tilde x,\tilde y)$.

\begin{lemma}
For any pairs $(\tilde x,\tilde y)$ and $(\tilde x',\tilde y')$ in $\R^2$, we have
$$\int_{M}\psi^k(\tilde x,\tilde y)dk=\int_{M}\psi^k(\tilde x',\tilde y')dk.$$
\end{lemma}
\proof If we consider the variable change $(x,y)\mapsto (x-\tilde x,y-\tilde y)$,
then
\begin{align*}
\int_{M}\psi^k(\tilde x,\tilde y)dk&=\int_{M}\frac{1}{\sigma^{3\slash2}}\psi_0(A^{-1}_{\theta,\sigma}(\tilde x-x,\tilde y-y))dxdyd\theta d\sigma\\
&= \int_{M}\frac{1}{\sigma^{3\slash2}}\psi_0(A^{-1}_{\theta,\sigma}(-x,-y))dxdyd\theta d\sigma\\
&=\int_{M}\psi^k(0,0)dk,
\end{align*}
and this prove the independence from $(\tilde x,\tilde y)$.\fine

In order to prove the convergence of (\ref{int}), we rewrite it by
making the change to polar coordinates $(\tilde x,\tilde y)\mapsto (r,\alpha)$,
such that 
$$\tilde x=r\cos(\al),\ \tilde y=r\sin(\al).$$
In order to have a simpler notation, we write $(\tilde x_\theta,\tilde y_\theta)=R_\theta(\tilde x,\tilde y)$,
that is $\tilde y_\theta=r\sin(\al+\theta)$.
In this notation we have
\begin{align*}
\int_{M}\psi^k(0,0)dk&=\int_{M}\frac{1}{\sigma^{3\slash2}}\cdot e^{-\sigma^{-2}(\tilde x^2+\tilde y^2)}\sin(-2\sigma^{-1}\tilde y_{\theta})
d\tilde xd\tilde y d\theta d\sigma\\
&= \int_{M}\frac{re^{-\sigma^{-2}r^2}}{\sigma^{3\slash 2}}\sin\left(-\frac{2r}{\sigma}\sin(\al+\theta)\right)drd\al d\theta d\sigma.
\end{align*}

\begin{prop}
The following integral converges to $0$,
$$\int_{M}\psi^k(0,0)dk.$$
\end{prop}
\proof We rewrite again the integral as
$$\int_{S^1\times \R^+}d\theta\sigma^{1\slash2}d\sigma \int_{S^1\times \R^+} \frac{1}{\sigma}drd\al \frac{re^{-\sigma^{-2}r^2}}{\sigma}
 \sin\left(-\frac{2r}{\sigma}\sin(\al+\theta)\right).$$
We consider the coordinate changes $\al'= \al-\theta$
and $s=\frac{r}{\sigma}$, obtaining the integral
\begin{equation}\label{eq_intint}
\int_{S^1\times \R^+}2\pi\sigma^{1\slash2}d\sigma\int_{S^1\times \R^+}dsd\al'se^{-s^2}
\sin\left(-2s\sin{\al'}\right).
\end{equation}
We have
$$\left|\int_{S^1\times \R^+}dsd\al' se^{-s^2}\sin\left(-2s\sin\al'\right)\right|\leq 2\pi\int_0^\infty se^{-s^2}ds=\pi,$$
therefore the internal integral in (\ref{eq_intint}) has finite value.
As the sinus is an odd function, its value is $0$.
Therefore by the Fubini's theorem
we have $\int_M\psi^k(0,0)=0$ too.\fine

\section{}\label{appb}
The goal of this appendix is to prove that the integral
\begin{equation}\label{intb}
C_{\psi}=\int_{\R^2}\frac{\left|\widehat \psi_0(\xi)\right|^2}{|\xi|^2}d\xi
\end{equation}
is finite and for any $f,g\in L^2(\R^2)$,
\begin{equation}\label{eqb}
\int_{M}\frac{dk}{\sigma^2}\langle f,\psi^k\rangle \langle \psi^k,g\rangle=C_\psi\cdot\langle f,g\rangle.
\end{equation}

Observe preliminarily that
$\widehat{\psi_0}(0)=0$. Indeed,
$$\int_{\R^2} e^{-\tilde x^2-\tilde y^2}\sin(2\tilde y)d\tilde xd\tilde y=0,$$
because the sinus is  an odd function.

\begin{lemma}
The integral (\ref{intb}) is finite.
\end{lemma}
\proof
We use the notation $\xi=(\xi_1,\xi_2)$.
Moreover, we observe $\sin(2\tilde y)=\frac{1}{2i}(e^{2i\tilde y}-e^{-2i\tilde y})$.
Therefore, we have
$$\widehat{\psi_0}(\xi)=\int_{\R^2}\psi_0(\tilde x,\tilde y)e^{-2\pi i (\tilde x,\tilde y)\cdot \xi}d\tilde x d\tilde y=
\int_{\R^2}e^{-\tilde x^2-\tilde y^2-2\pi i \xi_1\tilde x-2\pi i \xi_2\tilde y}\cdot\frac {e^{2i\tilde y}-e^{-2i\tilde y}}{2i}d\tilde xd \tilde y.$$
By completing the squares, we obtain
$$\widehat{\psi_0}(\xi)=\frac{1}{2i}e^{-\pi^2\xi_1^2}\int_{\R^2}e^{-(\tilde x+i\pi\xi_1)^2}\left(e^{-(\pi+1)^2\xi_2^2}e^{-(\tilde y+i(\pi +1)\xi_2)^2}
-e^{-(\pi-1)^2\xi_2^2}e^{-(\tilde y+i(\pi-1)\xi_2)^2}\right)d\tilde xd\tilde y.$$
For $c\in \R$, we denote by $D(c)$ the value
$$D(c):=\int_\R e^{-(\tilde x+ic)^2}d\tilde x,$$
which is known to be finite.
Therefore 
$$\widehat{\psi_0}(\xi)=\frac{1}{2i}e^{-\pi^2\xi_1^2}D(\pi\xi_1)\left(e^{-(\pi+1)^2\xi_2^2}D((\pi+1)\xi_2)-e^{-(\pi-1)^2\xi_2^2}D((\pi-1)\xi_2)\right).$$
We know that $D(0)=\sqrt\pi$, the integral of the Gaussian function.
Therefore the development of $|\widehat{\psi_0}|$ around $0$ is
$$\left|\widehat{\psi_0}(\xi)\right|=2\pi^2\xi_2^2+ o(|\xi|^3).$$
Moreover $D(c)\leq D(0)$ for any $c\in \R$, then
at infinity the function behaves as the difference of two Gaussians.
This proves that
$$\int_{\R^2}\frac{\left|\widehat{\psi_0}(\xi)\right|^2}{|\xi|^2}d\xi$$
is a convergent integral.\fine

\begin{prop}
The equality (\ref{eqb}) is satisfied for any pair of functions 
$f,g\in L^2(\R^2)$.
\end{prop}
\proof We start by observing that, for any $k=(x,y,\theta,\sigma)$,
\begin{equation}\label{ok3}
\widehat{\psi^k}(\xi)=\sigma^{1\slash2}\cdot\widehat{\psi_0}(\sigma R_{-\theta}\xi)\cdot e^{-2\pi i((x,y)\cdot \xi)}.
\end{equation}
Indeed, with the coordinate change $(\tilde x',\tilde y')=\sigma^{-1}R_{-\theta}(\tilde x-x,\tilde y-y)$ we have
\begin{align*}
\widehat{\psi^k}(\xi)&=\int_{\R^2} \psi_0(\sigma^{-1}R_{-\theta}(\tilde x-x,\tilde y-y))e^{-2\pi i((\tilde x,\tilde y)\cdot \xi)}d\tilde xd\tilde y\\
&=\left(\int_{\R^2}\psi_0(\tilde x',\tilde y')e^{-2\pi i((\tilde x',\tilde y')\cdot (\sigma R_{-\theta}\xi))}\sigma^{1\slash2}d\tilde x' d\tilde y'\right)e^{-2\pi i ((x,y)\cdot\xi)}d\tilde xd\tilde y,\\
\end{align*}
and therefore we obtain (\ref{ok3}). 

Moreover we observe that
$$\int_{\R^2}\widehat{\psi^k}(\xi)\cdot\overline{ f(\xi)}d\xi=\int_{\R^2}\sigma^{1\slash2}\cdot\widehat{\psi_0}(\sigma R_{-\theta}\xi)\cdot \overline{ f(\xi)}
\cdot e^{-2\pi i((x,y)\cdot \xi)}d\xi $$
is the Fourier transform of $F_f(\xi):=\sigma^{1\slash2}\cdot\widehat{\psi_0}(\sigma R_{-\theta}\xi)\cdot\overline{ f(\xi)}$ evaluated 
in $(x,y)$.

In order to
prove (\ref{eqb}), we observe that by the Plancherel Theorem,
\begin{align*}
\int_{M}\langle f,\psi^k\rangle\langle \psi^k,g\rangle\frac{dk}{\sigma^2}&=
\int_{M}\frac{d\theta d\sigma}{\sigma^2}\left(\int f(\xi)\overline{\widehat{\psi^k}(\xi)}d\xi\right)\left(\int \widehat{\psi^k}(\xi)\overline {g(\xi)}d \xi\right)\\
&=\int_{M}\frac{d\theta d\sigma}{\sigma^2}\overline{\widehat{F_f}(x,y)}\widehat{F_g}(x,y)d\xi\\
&=\int_{M}\frac{d\theta d\sigma}{\sigma^2}\sigma\cdot f(\xi)\cdot\overline{g(\xi)}\cdot\left|\widehat{\psi_0}(\sigma R_{-\theta}\xi)\right|^2d\xi.
\end{align*}
By Fubini's Theorem we have
\begin{align*}
\int_{M}\langle f,\psi^k\rangle\langle \psi^k,g\rangle\frac{dk}{\sigma^2}&=
\int_{\R^2}d\xi f(\xi)g(\xi)\int_0^{2\pi}\int_0^\infty\frac{d\theta d\sigma}{\sigma}\left|\widehat{\psi_0}(\sigma R_{-\theta}\xi)\right|^2\\
&=\int_{\R^2}d\xi f(\xi)g(\xi)\int_{\R^2}\frac{\left|\widehat{\psi_0}(h)\right|^2}{|h|}dh\\
&=\langle f,g\rangle\cdot C_\psi,
\end{align*}
where we used the change of coordinates 
$(\theta,\sigma)\mapsto h=\sigma R_{-\theta}\xi$.\fine

\section{}\label{appc}

In this appendix we prove that equality (\ref{eq_eq0})
holds not only in the weak sense, but also in a much stronger version.
In particular, for any $\sigma_1,\sigma_2, B\in \R^+$ and $f\in L^2(\R^2)$, we define
$$\int_0^{2\pi}d\theta\int_{\sigma_1}^{\sigma_2}\frac{d\sigma}{\sigma^2}\int_{\norm{(x,y)}\leq B}dx dy \langle f, \psi^k\rangle \psi^k$$
as the unique element in $L^2(\R^2)$ whose inner product with any $g\in L^2(\R^2)$ is
$$\int_0^{2\pi}d\theta\int_{\sigma_1}^{\sigma_2}\frac{d\sigma}{\sigma^2}\int_{\norm{(x,y)}\leq B}dx dy \langle f, \psi^k\rangle \langle\psi^k, g\rangle.$$
\begin{prop}
For any $f\in L^2(\R^2)$,
$$\lim_{\begin{array}{c}
\sigma_1\to 0\\ 
\sigma_2,B\to \infty\\
\end{array}}\norm{f-C_\psi^{-1}\int_0^{2\pi}\int_{\sigma_1}^{\sigma_2}\int_{\norm{(x,y)}\leq B}
\langle f,\psi^k\rangle \psi^k \frac{1}{\sigma^2}dxdyd\theta d\sigma}=0.$$
\end{prop}
\proof The term we are estimating is
$$ \sup_{\norm{g}=1}\left| \left\langle f-C_\psi^{-1}\int_0^{2\pi}\int_{\sigma_1}^{\sigma_2}\int_{\norm{(x,y)}\leq B}
\langle f,\psi^k\rangle \psi^k\frac{dk}{\sigma^2}, \ g\right\rangle   \right|$$
which is bounded above by
$$ \sup_{\norm{g}=1}\left| C_\psi^{-1}\int_{k\in U_{\sigma_1,\sigma_2,B}}\langle f,\psi^k\rangle \langle \psi^k,g\rangle \frac{dk}{\sigma^2}\right|,$$
where $U_{\sigma_1,\sigma_2, B}\subset M$ is defined by
$$U_{\sigma_1,\sigma_2,B}=\left\{k\in M|\ \sigma<\sigma_1\m{ or }\sigma>\sigma_2\m{ or }\norm{(x,y)}> B\right\}.$$
By Caucy-Schwarz inequality this is again bounded by
$$\leq \sup_{\norm{g}=1}\left|C_\psi^{-1}\int_{k\in U_{\sigma_1,\sigma_2,B}}|\langle f,\psi^k\rangle|^2\frac{dk}{\sigma^2}\right|^{1\slash 2}\cdot
\left| C_\psi^{-1}\int_M|\langle g,\psi^k\rangle |^2\frac{dk}{\sigma^2}\right|^{1\slash 2}.$$
The second term is $\norm{g}=1$ by equality (\ref{eqb}).
The first term tends to $0$ for $\sigma_1\to 0$ and $\sigma_2,B\to \infty$,
because the infinite integral $\int_M|\langle f,\psi^k\rangle|^2\frac{dk}{\sigma^2}$ converges to $\norm{f}^2$.\fine

\section{}\label{appd}
In order to evaluate the distance $d_c$ between two points $p_0=(x_0,y_0,\theta_0,\sigma_0)$
and $p_1=(x_1,y_1,\theta_1,\sigma_1)$, we search for the constant
coefficients that allows a vector field $Y$ whose flow relies $p_0$ to $p_1$.

We start by observing that 
$$\dot\theta=c_2,$$
therefore $c_2=\theta_1-\theta_0$ and
$$\theta_t=c_2\cdot t+\theta_0.$$
Regarding the $\sigma$ variable, we have
$$\dot\sigma=\frac{c_4}{h_2}\cdot \sigma,$$
therefore $c_4=\ln\left(\frac{\sigma_1}{\sigma_0}\right)\cdot h_2$
and
$$\sigma_t=\sigma_0\cdot e^{\frac{c_4t}{h_2}}.$$
In what follows we will use the notation $\tilde c_4:=\frac{c_4}{h_2}$.

Regarding the first two coordinates
$$\binom{\dot x}{\dot y}=\left(
\begin{array}{cc}
\cos\theta &-\frac{\sigma}{h_1}\cdot \sin\theta\\
\sin\theta & \frac{\sigma}{h_1}\cdot \cos\theta
\end{array}\right)\binom{c_1}{c_3}.$$

We denote by $S_t$ the matrix obtained by integrating the matrix 
above up to time $t$, therefore
$$S_t=\left(\begin{array}{cc}
S_t^{11} & S_t^{12}\\
S_t^{21} & S_t^{22}
\end{array}\right)$$
with
\begin{align*}
S_t^{11}&=\left[\frac{\sin(\theta_t)}{c_2}\right]^t_0\\
S_t^{12}&=\left[\frac{1}{h_1\cdot(\tilde c_4^2+c_2^2)}\cdot \left(c_2\cdot \sigma_t\cos(\theta_t)
-\tilde c_4\cdot \sigma_t\sin(\theta_t)\right)\right]^t_0\\
S_t^{21}&=\left[\frac{-\cos(\theta_t)}{c_2}\right]^t_0\\
S_t^{22}&=\left[\frac{1}{h_1\cdot (\tilde c_4^2+c_2^2)}\cdot \left(\tilde c_4\cdot \sigma_t\cos(\theta_t)+c_2\cdot 
\sigma_t\sin(\theta_t)\right)\right]^t_0.
\end{align*}\newline

We obtain 
$$\binom{c_1}{c_3}=S_1^{-1}\binom{x_1-x_0}{y_1-y_0}$$
and
$$\binom{x_t}{y_t}=S_t\binom{c_1}{c_3}+\binom{x_0}{y_0}.$$

For any two points as above, their  distance is
$$d_c(p_1,p_0)^2=c_1^2+c_2^2+c_3^2+c_4^2.$$\newline

\cor{The datasets and codes generated during and analysed during the current study are available from the corresponding author on reasonable request.}

\bibliography{mybiblio}{}
\bibliographystyle{plain}

\end{document}